\documentclass[a4paper,11pt]{article}
\usepackage{amsfonts,amssymb,latexsym}
\begin{document}

\Large

\newcommand{\Dp}{\Delta^+}
\newcommand{\De}{\Delta}

\newcommand{\Zb}{\mathbb{Z}}
\newcommand{\Cb}{\mathbb{C}}
\newcommand{\Pb}{\mathbb{P}}

\newcommand{\gog}{{\mathfrak g}}
\newcommand{\nog}{{\mathfrak n}}
\newcommand{\mog}{{\mathfrak m}}
\newcommand{\hog}{{\mathfrak h}}
\newcommand{\bog}{{\mathfrak b}}
\newcommand{\pog}{{\mathfrak p}}
\newcommand{\ut}{{\mathfrak u}{\mathfrak t}}
\newcommand{\gl}{{\mathfrak g}{\mathfrak l}}
\newcommand{\sg}{\sigma}
\newcommand{\eps}{\varepsilon}
\newcommand{\ueps}{\epsilon}
\newcommand{\ad}{{\mathrm{ad}}}
\newcommand{\id}{{\mathrm{id}}}
\newcommand{\IC}{{\cal I}}
\newcommand{\MC}{{\cal M}}
\newcommand{\CC}{{\cal C}}
\newcommand{\SC}{{\cal S}}

\newcommand{\ii}{{\mathfrak i}}
\newcommand{\jj}{{\mathfrak j}}
\newcommand{\DC}{{\cal D}}
\newcommand{\AC}{{\cal A}}
\newcommand{\JC}{{\cal J}}
\newcommand{\GL}{\mathrm{GL}}
\newcommand{\UT}{\mathrm{UT}}
\newcommand{\Ad}{{\mathrm{Ad}}}
\newcommand{\Ann}{{\mathrm{Ann}}}
\newcommand{\ord}{{\mathrm{ord}}}
\newcommand{\lee}{{\leqslant}}
\newcommand{\gee}{{\geqslant}}
\newcommand{\La}{{\Lambda}}

\title{ Involution in   $S_n$  and associated coadjoint orbits}
\date{}
\author{A.N.Panov\\ apanov@list.ru\thanks{The paper is supported by  grants RFBR
  05-01-00313  and 06-01-00037}}
\maketitle
\section{Introduction and main statements }

 Coadjoint orbits  play an important role in
representation theory, harmonic analysis and theory of symplectic
varieties ~\cite{K-Orb,D}. For semisimple Lie groups coadjoint
orbits coincide with  with adjoint orbits and their classification
is  well known.  The classical example of  nilpotent Lie group is
the unitriangular group. The problem of classification of coadjoint
orbits for this group is still far from its solution
~\cite{K-Orb,K-Var}. Description of regular orbits (i.e. orbits of
maximal dimension) for $\UT(n,K)$ was achieved by A.A.Kirillov in
the his origin paper on the orbit method \cite{K-62}. Subregular
orbits for an arbitrary $n$ and all orbits for  $n\lee 7$ are
classified in ~\cite{IP}.

The goal of this paper is to enlarge the class of orbits  that have
 complete description. To each involution $\sg$ in
the symmetric group $S_n$ we correspond a family of coadjoint orbits
$\{\Omega(f)\vert~ f\in X_\sg\}$. We construct a polarization of the
canonical form of any such orbit (Theorem 1.1). We have got the
formula for dimension of these orbits (Theorem 1.2).

Orbits of a  regular action of an algebraic unipotent group on an
affine  algebraic variety are closed ~\cite[11.2.4]{D}. In the paper
we find generators in the defining ideal of the  orbit
$\Omega(f)$(Theorem 1.4).
 Note that each of constructed generators has  the form  $P-c$, where
  $P$ is a coefficient of minor of characteristic matrix  and $c\in K$.
 The class of orbit associated with involutions contains all regular and some
 subregular orbits( see Examples 2,3 in \S1 and ~\cite{IP}).
 Remark also that one can set up a   correspondence  between    involutions  and families
 of  adjoint orbits   satisfying  the equation  $X^2=0$ ~\cite{Mel}.

Now go on to formulation of the main statements of this paper. An
unitriangular group is a group  $N:=\UT(n,K)$ which consists of  all
lower triangular matrices over a field $K$ with units on the
diagonal. We suppose that the field $K$ has zero characteristic.
   The
 Lie algebra  $\nog=\ut(n,K)$  of group $N$ consists of the lower
 triangular matrices with zeros on the diagonal.
 With the help of Killing form $(\cdot,\cdot)$
 we identify  the conjugate space $\nog^*$ with the subspace
 $\nog_+$ that consists of the upper triangular matrices with zeros on the  diagonal.
 We also identify the symmetric algebra  $S(\nog)$ with the algebra
 $K[\nog^*]$ of regular functions on
 $\nog^*$.
The group $N$ acts on  $\nog^*$ by the formula $\Ad^*_gf(x) =
f(\Ad^{-1}_gx)$,~ $g\in N, f\in\nog^*$.

A subalgebra  $\pog$ in $\nog$ is a polarization of  $f\in \nog^*$
if  $\pog$ is maximal isotropic subspace with respect to the skew
symmetric form $f([x,y])$. According to the orbit method  for  any
linear form on a nilpotent Lie algebra there exists a polarization.
Existence  of polarization gives the possibility to construct
 a primitive ideal in the universal enveloping algebra
 $U(\nog)$ ~\cite{D} and in the case  $K=\mathbb{R}$
the irreducible representation of the group $N$~\cite{K-Orb}.

Denote by  $\Dp$ the system of positive roots of the Lie algebra
$\mathrm{gl}(n,K)$ ~\cite{GG}. the set of positive roots. The
relation $>$ in $\Dp$ is  a lexicographical order. Any positive root
$\xi$ has the form $\alpha_{ji} = \eps_j-\eps_i$, where $j<i$. We
denote $j=j(\xi)$ and $i=i(\xi)$. Any positive $\xi$ root defines
the reflection $r_\xi\in S_n$.

 Let $\sigma$ be the
involution in  $S_n$ (i.e. $\sigma^2=\mathrm{id}$). The involution
$\sigma$ uniquely decomposes into a product of commutating
reflections
$$\sigma = r_1r_2\cdots r_s,\eqno(1)$$ where each $r_m$ is the reflection
with respect to the positive  root $\xi_m$. We assume that $\xi_1>
\xi_2>\ldots>\xi_s$. Since reflections mutually commute, then
$j(\xi_1)<j(\xi_2)<\ldots<j(\xi_s)$. Form the subset
$$\SC=\{\xi_1, \xi_2, \ldots,\xi_s\}.$$
Let  $\{y_{ij}\}_{1\lee j < i\lee n-1}$ be the standard basis
$\nog$. In what  follows we shall also  use the notation $y_\gamma$
for  $y_{ij}$  if $\gamma=\eps_j-\eps_i$.

 Consider the subset  $X_\sigma\subset \nog^*$,
which  consists of   $f\in{\mathfrak n}^*$ that $f(y_{\xi_m})\ne 0$,
$1\lee m\lee s$, and   $f$ annihilates  on the other  vectors of
standard basis. \\
{\bf Main definition}. We call the orbit $\Omega(f)$ of   $f\in
X_\sg$ a coadjoint orbit associated with the involution $\sg$.

Our first goal is to construct a polarization of  any $f\in X_\sg$.
One can decompose the set of positive roots:
 $$\Dp=\bigsqcup_{1\lee t\lee n-1} \De^{(t)},\quad \mbox{where}~~
\De^{(t)}=\{ \gamma\in\Dp\vert~ j(\gamma)=t\}.$$ For any $1\lee t
\lee n-1$ we denote
$$\sg_t =\prod_{j(\xi_m)\lee t} r_m.$$
Easy to see that  $\sg_{n-1 } = \sg$. We put  $\sg_0=\id$. Denote by
$\Pi_\sg$ the following subset of positive roots

$$\Pi:=\Pi_\sg:= \bigsqcup_{1\lee t\lee n-1} \Pi^{(t)},$$
$$ \mbox{where}~~ \Pi^{(t)}=\{ \eta\in\De^{(t)}\vert~ \sg_{t-1}(\eta) > 0\}.$$
Denote by $\pog_\sg$  the  linear subspace spanned by $\{ y_{it}:~\eps_t-\eps_i\in\Pi\}.$\\
 {\bf Theorem 1.1}.  $\pog_\sg$ is a polarization for any $f\in X_\sg$.

 As usual  $l(\sg)$ (resp. $s(\sg)$)  is a number of simple (resp. arbitrary) reflections in reduced decomposition of
  $\sg$  into a product of  simple(resp. arbitrary) reflections.
  Easy to see that  $s(\sg) = s$ (see (1)).\\
 {\bf Theorem 1.2}.  For any  $f\in X_\sg$ the
 dimension of the orbit  $\Omega(f)$  equals to  $l(\sg) - s(\sg)$.

Our next goal  is to find  system of generators in the defining
ideal $\IC(\Omega(f))$ of the orbit  $\Omega(f)$  (i.e. in the ideal
that   consists of all polynomials that annihilate on
 $\Omega(f)$).
Consider the following formal matrix  $\Phi$ and its characteristic matrix $\Phi(\tau)$:

$$\Phi = \left(\begin{array}{cccc}0&0&\ldots&0\\
y_{21}&0&\ldots&0\\
\vdots&\vdots&\ddots&\vdots\\
y_{n1}&y_{n2}&\ldots&0
\end{array}
\right),\quad \Phi(\tau) = \tau\Phi+E =
\left(\begin{array}{cccc}1&0&\ldots&0\\
\tau y_{21}&1&\ldots&0\\
\vdots&\vdots&\ddots&\vdots\\
\tau y_{n1}&\tau y_{n2}&\ldots&1
\end{array}
\right).
$$
For any subsets $J,I\subset\{1,\ldots,n\}$,~  $|I| = |J|$, we denote
by $M_I^J$ the minor of  matrix $\Phi$ with ordered increasing
systems of columns $\mathrm{ord}(J)$ and rows $\mathrm{ord}(I)$.

 We denote by  $M^J_I(\tau)$ the
corresponding minor of matrix $\Phi(\tau)$. The minor  $M_I^J(\tau)$
is a polynomial in   $\tau$ with coefficients in the symmetric
algebra $S(\nog)$:
$$M_I^J(\tau) =
\tau^m(P_{I,0}^J + P_{I,1}^J\tau+ P_{I,2}^J\tau^2+ \ldots
),\eqno(2)$$ where  $m$  equals to $\vert J\setminus I\vert = \vert
I\setminus J\vert$. The   coefficient  $P_{I,0}^J$ is up to a sign
equals to
 the minor $M^{J\setminus I}_{I\setminus
J}$ of the matrix $\Phi$ (in what follows we frequently drop this
sign).
\\
{\bf Remark 1}. If $I\gee J$ (in the sense Definition 3.4), then
$M_I^J(\tau)\ne 0$ and its first coefficient $P_{I,0}^J$ is also
nonzero. If $I\ngeqslant J$, then $M_I^J(\tau)= 0$.

Note that  for any  $f\in \nog_+=\nog^*$ the value of  $M_I^J$ at
the point $f$ coincides with the transposed minor of $f$ presented
as an upper triangular matrix. Similarly, the values of $P_{I,0}^J,
P_{I,1}^J, P_{I,2}^J, \cdots $  at the point $f$ coincides with
values of coefficients of the corresponding minor of the
characteristic matrix $\tau f+E$.

For any pair   $1\le k,t\le n$ we consider the  systems
$$ J'(k,t)=\{1\le j < t: ~ \sigma(j) > k\},\quad\quad
 I'(k,t)=\sigma J'(k,t),$$
$$J(k,t)=
J'(k,t)\sqcup\{t\},\quad\quad I(k,t)= I'(k,t) \sqcup \{k\}.$$

Denote by  $D_{k,t}$ (resp. $D_{k,t}(\tau)$) the  minor of the
matrix $\Phi$ (resp. $\Phi(\tau)$) with the systems of columns
$J(k,t)$ and rows $I(k,t)$.

As in (2) we decompose $$D_{k,t}(\tau)= \tau^m( P_{k,t,0}+
P_{k,t,1}\tau + P_{k,t,2} \tau^2 + \ldots),\eqno(3)$$ where the
coefficients $P_{k,t,i}$ are elements of $ S({\mathfrak
n})=K[{\mathfrak n}^*]$ and $P_{k,t,0}\ne 0$. Note that  if  $k>t$
then $D_{k,t}(\tau)= \tau^m P_{k,t,0}$. Denote
$$ \MC =\{ \gamma\in \Dp\vert~
\sg_{t}(\gamma)>0,~\mbox{where}~ t=j(\gamma)\}.$$

To each root of $\MC$ we shall assign a polynomial on $\nog^*$.\\
 {\bf Definition 1.3}. Let $i>t$ and $\eta=\eps_t-\eps_i\in\MC$.
We say that $\eta$ has type 0, if there is no $\xi\in\SC$ such that
$j(\xi)\lee t-1$ and $i(\xi)=i$. Otherwise, we say that a $\eta$ of
$\MC$ has type 1.

Put $k=\sg_{t-1}(i)$.
 We denote
$$Q_{i,t}= \left\{\begin{array}{l} P_{k,t,0},~\mbox{if}~ \eta~\mathrm{has~type~0},
\\
P_{k,t,1},~\mbox{if}~ \eta~\mathrm{has~type~1}.
\end{array}\right.$$

Note that, if $\eta$  has type 0, then $k=i$ and $Q_{it}=P_{i,t,0}$
coincides (up to a sign) with the minor $D_{i,t}$ of matrix $\Phi$.
If $\eta$ has type 1, then $k<t$  and $k=\sg(t)$. For each $1\lee
m\lee s$ we denote
$$D_m =
 D_{i(\xi_m),j(\xi_m)}.$$
 Easy to see that    $D_m(f)\ne 0$ for any  $f\in X_\sg$.\\
{\bf Theorem 1.4}. For any $f\in X_\sigma$ the defining ideal of
coadjoint orbit
 $\Omega(f)$
is generated by  $Q_{i,t}$, where $\eps_t-\eps_i\in \MC$, and $D_m -
D_m(f)$, where   $ 1\lee m\lee s$. \\
{\bf Corollary 1.5}. For any  $f\in X_\sg$ the intersection
$\Omega(f)\cap X_\sg$ contains the only element $f$.

 For each $1\lee t\lee n-1$ decompose the involution  $\sg=\sg_t\sg_t'$ where
$$\sg_t'=\prod_{j(\xi_m)>t} r_m.$$
Denote $$\Dp_\sg =\{\zeta\in\Dp:~ \sg(\eta)>0\}.$$ By definition,
$|\Dp_\sg|=|\Dp|-l(\sg)$.
\\
 {\bf Remark 2}. The positive root  $\zeta$ lies in $\De_\sg^{(t)}=\Dp_\sg\cap\De^{(t)}$
  if and only if
$\sg_t(\eta)>0$ where $\eta=\sg_t'(\zeta)$. This shows that
  $\sg_t' $  bijectively maps $\De_\sg^{(t)}$ onto $\MC^{(t)}=\MC\cap
 \De^{(t)}$. This implies $ |\MC^{(t)}| = |\De_\sg^{(t)}|$
 and, therefore, $ |\MC| = |\Dp_\sg|= |\Dp|-l(\sg)$.
  From Remark 2 and Theorems 1.2 and 1.4 we obtain\\
 {\bf Corollary 1.6}. Each  coadjoint orbit $\Omega(f)$, $f\in
 X_\sg$ is a complete intersection.

 Theorems 1.1 and 1.2 are proved in  \S 2.
The next \S 3 is completely devoted to the proof of Theorem 1.4. In
the sequel of this section we introduce the language of admissible
diagrams and reformulate Theorems 1.1 in new terms. This approach
has its own advantages. We present examples of involutions,
construct
  their admissible diagrams and generators of defining ideals of their associated orbits.

As above $\sg$ be an involution in  $S_n$. Consider the following
subsets in the set of positive roots $\Dp$:
$$ \CC_+=\{ \gamma\in \Dp\setminus \SC\vert~ \sg_{t-1}(\gamma)>0,
\sg_{t}(\gamma)<0, ~\mbox{where}~t=j(\gamma)\},$$
$$\CC_-=\{\gamma\in\Dp\vert ~ \sg_{t-1}(\gamma)<0, ~\mbox{where}~t=j(\gamma)\}$$
The set of positive roots  $\Dp$ and its subset $\Pi$ decomposes:
$\Dp= \Pi\sqcup \CC_-$  and $\Pi=\SC\sqcup \MC\sqcup \CC_+$.

 We  begin construction of an admissible diagram from the empty $n\times n
$-matrix. We shall not fill the  squares  $(i,j)$, $i \lee j$. To
each positive root $\gamma=\alpha_{ji}$, $i>j$, we correspond the
square $(i,j)$.  We put one of the symbols $\otimes$, $+$, $-$,
$\bullet$ into  square $(i,j)$  according to these rules:
$$\left\{\begin{array}{l} \otimes ~, ~~\mbox{if}~ \gamma\in \SC,\\
 + ~,~~\mbox{if}~ \gamma\in \CC_+,\\
 - ~,~~\mbox{if}~ \gamma\in \CC_-,\\
  \bullet ~, ~~\mbox{if}~ \gamma\in M.\\
 \end{array}\right.$$

One can construct the same  admissible diagram using the method of
paper $~\cite{IP}$. Consider the first root  $\xi_1$  in  $\SC$. Put
the symbol $\otimes$ into  square  $(i(\xi_1), j(\xi_1))$. Next we
put the symbol $+$ into all  squares that lie in the column
$j(\xi_1)$ lower the diagonal and higher than the  square
$(i(\xi_1), j(\xi_1))$. We put the symbol $-$ into all  squares that
lie in the row $i(\xi_1)$ lower the diagonal and righter than the
square $(i(\xi_1), j(\xi_1))$.

After this procedure,  we take the other root $\xi_2$ of $\SC$.   We
put the symbol $\otimes$ into   the  square $(i(\xi_2), j(\xi_2))$
that as it easy to see are not filled after the first step. Similar
to the case of root $\xi_1$ we put the symbols $+$ and $-$ into the
matrix using the following additional rule: if for certain $k$ one
of this two squares $(k,j(\xi_2))$, $(i(\xi_2),k)$ are already
filled we do not fill the other  square.

We continue this procedure for all other roots  $\xi\in \SC$. We put
the symbol $\bullet$ into all unfilled  squares lower the diagonal.
Note that  the roots of  $\Pi$ correspond to the squares filled by
symbols $\otimes$, $+$ and $\bullet$.
\\
{\bf Example 1(begin)}. Involution  $\sg=(1,4)(2,7)(3,6)$  in the
symmetric group $S_7$ corresponds to the admissible diagram

\begin{center} {{\large
\begin{tabular}{|p{0.1cm}|p{0.1cm}|p{0.1cm}|p{0.1cm}|p{0.1cm}|p{0.1 cm}|p{0.1 cm}|}
\hline & & & & & &\\
\hline $+ $& & & & & &\\
\hline $+$ & $+$ & & & & &\\
\hline $\otimes$ &$-$  & $-$ & & &&\\
\hline $\bullet$ & $+$ & $+$ & $\bullet$ & & &\\
\hline $\bullet$ & $+$ & $\otimes$ & $\bullet$ & $-$& & \\
\hline $\bullet$ & $\otimes$ & $-$ & $\bullet$ & $-$& $-$& \\
\hline
\end{tabular}}}
\end{center}
The diagram is constructed in 4 steps:
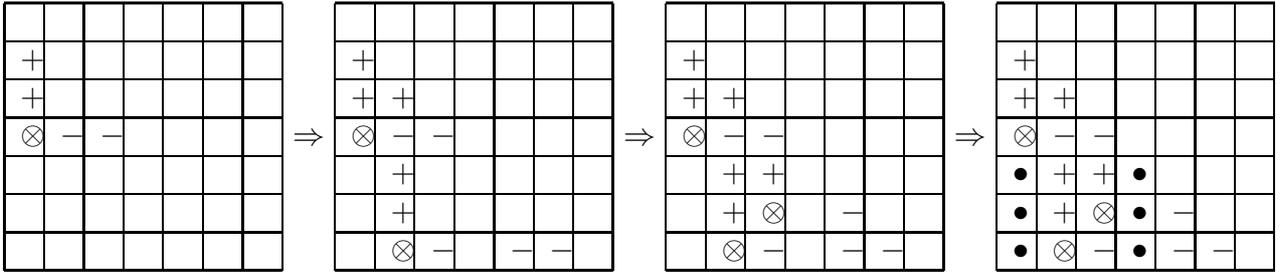
\begin{center} {{\large
\begin{tabular}{|p{0.1cm}|p{0.1cm}|p{0.1cm}|p{0.1cm}|p{0.1cm}|p{0.1 cm}|p{0.1 cm}|}
\hline & & & & & &\\
\hline $+ $& & & & & &\\
\hline $+$ &  & & & & &\\
\hline $\otimes$ &$-$  & $-$ & & &&\\
\hline  &  &  &  & & &\\
\hline  &  & &  & & & \\
\hline  &  &  &  & & & \\
\hline
\end{tabular}
$\Rightarrow$
\begin{tabular}{|p{0.1cm}|p{0.1cm}|p{0.1cm}|p{0.1cm}|p{0.1cm}|p{0.1 cm}|p{0.1 cm}|}
\hline & & & & & &\\
\hline $+ $& & & & & &\\
\hline $+$ & $+$ & & & & &\\
\hline $\otimes$ &$-$  & $-$ & & &&\\
\hline  & $+$ &  & & & &\\
\hline  & $+$ &  & & & & \\
\hline  &$\otimes$ & $-$ &  & $-$& $-$& \\
\hline
\end{tabular}
$\Rightarrow$
\begin{tabular}{|p{0.1cm}|p{0.1cm}|p{0.1cm}|p{0.1cm}|p{0.1cm}|p{0.1 cm}|p{0.1 cm}|}
\hline & & & & & &\\
\hline $+ $& & & & & &\\
\hline $+$ & $+$ & & & & &\\
\hline $\otimes$ &$-$  & $-$ & & &&\\
\hline  & $+$ & $+$ &  & & &\\
\hline  & $+$ & $\otimes$ &  & $-$& & \\
\hline  & $\otimes$ & $-$ &  & $-$& $-$& \\
\hline
\end{tabular}
$\Rightarrow$
\begin{tabular}{|p{0.1cm}|p{0.1cm}|p{0.1cm}|p{0.1cm}|p{0.1cm}|p{0.1 cm}|p{0.1 cm}|}
\hline & & & & & &\\
\hline $+ $& & & & & &\\
\hline $+$ & $+$ & & & & &\\
\hline $\otimes$ &$-$  & $-$ & & &&\\
\hline $\bullet$ & $+$ & $+$ & $\bullet$ & & &\\
\hline $\bullet$ & $+$ & $\otimes$ & $\bullet$ & $-$& & \\
\hline $\bullet$ & $\otimes$ & $-$ & $\bullet$ & $-$& $-$& \\
\hline
\end{tabular}
 }}
\end{center}
In this example  $\SC = \{\alpha_{14}, \alpha_{27}, \alpha_{36}\}$,
$\MC  = \{\alpha_{15}, \alpha_{16}, \alpha_{17}, \alpha_{45},
\alpha_{46}, \alpha_{47}\}$, $\Pi = \Dp\setminus \{
\alpha_{24},\alpha_{34},\alpha_{37},\alpha_{56}, \alpha_{57},
\alpha_{67}\}$. In the following Remarks 3-5 we reformulate Theorems
of this section in terms of  admissible diagrams.
\\
{\bf Remark 3}. According to Theorem 1.1.   the polarization
$\pog_\sg$ is spanned by  $\{ y_{it}\}$ such that the   square
$(i,t)$ is filled by one of the symbols  $\otimes$, $\bullet$ or $+$
in the
corresponding diagram. \\
{\bf Remark 4}. By Theorem 1.2 and Remark 2, $\dim\Omega(f)=l(\sg)
-s = |\Dp| - |M| - s$. We see that $\dim \Omega(f)$ is equal to the
number of $\pm$'s
in the admissible diagram.\\
{\bf Remark 5}. Let the  square $(i,t)$, $t<i$, be filled by
$\bullet$ (i.e. $\eps_t-\eps_i\in\MC$).
  We say that  $(i,t)$  has type 0(resp.  type 1), if
the root $\eps_t-\eps_i$ of $\MC$ has type 0(resp.  type 1). Easy to
see that $(i,t)$ has  type 0, if there is no symbol $\otimes$ in the
$i$th row and  on the left side of $(i,t)$. If the above symbol
$\otimes$ exists, then $(i,t)$ has type 1. In the Example 1, $(5,4)$
is type 0, but $(6,4)$  and $(7,4)$ have  type 1.

Note that, if  $(i,t)$ is filled by  symbol $\bullet$ and has  type
1, then $a < k < t <i$ where  $a=\sg(t)$ and as above  $k
=\sg_{t-1}(i)=\sg(i)$. The squares  $(t,a)$ and $(i,k))$ are filled
by symbol   $\otimes$. In Example 1, the square  $(7,4)$ is filled
by symbol  $\bullet$ and has  type 1, $a=1$ and  $k=2$. The squares
 $(4,1)$ and  $ (7,2)$ are filled by symbol $\otimes$. End of
 Remark.\\
{\bf Example 1(end)}. The defining  ideal  $\IC(\Omega(f))$  of
Example 1 is generated by  the elements, $D_m-c_m$, where $1\lee
m\lee 3$ and $c_m\ne 0$, and $Q_{i,t}$ corresponding to pairs
$\{(i,t)\}$ filled by symbol  $\bullet$. Here are the generating
polynomials of Example 1:
$$
D_1 = y_{41},\quad Q_{5,1} = y_{51},\quad Q_{6,1} = y_{61},\quad
Q_{7,1} = y_{71},$$
$$
D_2 = y_{72}, \quad
D_3= \left\vert\begin{array}{cc} y_{62}&y_{63}\\
y_{72}&y_{73}\end{array}\right\vert, \quad
Q_{5,4} = \left\vert\begin{array}{ccc} y_{52}&y_{53}&y_{54}\\
y_{62}&y_{63}&y_{64}\\
y_{72}&y_{73}&y_{74}
\end{array}\right\vert,$$
$$
Q_{6,4}= \left\vert\begin{array}{cc} y_{62}&y_{64}\\
y_{72}&y_{74}\end{array}\right\vert y_{41} +
\left\vert\begin{array}{cc} y_{62}&y_{63}\\
y_{72}&y_{73}\end{array}\right\vert y_{31}, \quad Q_{7,4} = y_{74}
y_{41} + y_{73} y_{31} + y_{72} y_{21}.
$$
{\bf Example 2}.  $\sg$  is an substitution of maximal length, i.e
$\sg=(1,n)(2,n-1)\cdots (k,k+1)$ where
 $k=[\frac{n}{2}]$. The
associated coadjoint orbit is a regular orbit (i.e. orbit of maximal
dimension). The defining ideal $\IC(\Omega(f))$ is generated by
$D_{i,n-i+1} - D_{i,n-i+1}(f)$, $1\lee i\lee k$.\\
{\bf Example 3}.  $\sg =(1,n-1)(2,n)(3,n-3)\cdots (k,k+1)$. The
associated coadjoint orbit of this involution is a subregular orbit
(i.e. orbit of dimension $\dim(\mathrm{reg.orbit}) -2$)
(see~\cite{IP}). The defining ideal $\IC(\Omega(f))$ is generated by
$$D_{i,n-i+1} - D_{i,n-i+1}(f),~3\lee i\lee k,$$
$$y_{n1},~~ y_{n-1,1}-f(y_{n-1,1}),~~y_{n2}- f(y_{n2}),~~z
=y_{n,n-1}y_{n-1,1}+y_{n,n-2}y_{n-2,1}+\ldots+y_{n2}y_{21}.$$ Here
are  the admissible diagrams for Examples 2 and 3 for $n=6$:

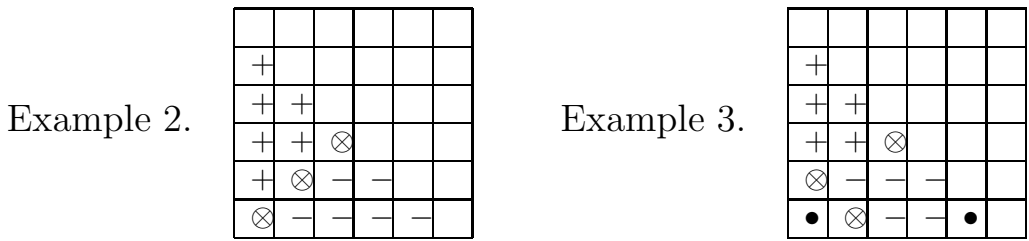
\begin{center} { Example 2.~~
{\large
\begin{tabular}{|p{0.1cm}|p{0.1cm}|p{0.1cm}|p{0.1cm}|p{0.1
cm}|p{0.1 cm}|}
\hline  & & & & &\\
\hline  +& & & & &\\
\hline   +&+ & & & &\\
\hline  + &+  & $\otimes$ & & &\\
\hline   + & $\otimes$  &$-$ &$-$ & &\\
\hline  $\otimes$ &$-$&$-$  & $-$&$-$ & \\
 \hline
\end{tabular}}
\quad\quad Example 3.~~ {\large
\begin{tabular}{|p{0.1cm}|p{0.1cm}|p{0.1cm}|p{0.1cm}|p{0.1 cm}|p{0.1 cm}|}
\hline  & & & & &\\
\hline  +& & & & &\\
\hline   +&+ & & & &\\
\hline  + &+  & $\otimes$ & & &\\
\hline   $\otimes$ & $-$  &$-$ &$-$ & &\\
\hline  $\bullet$ &$\otimes$&$-$  & $-$&$\bullet$ & \\
 \hline
\end{tabular}}
}\end{center}

\section{Polarization and dimension of orbit  }

The goal of this section is to prove Theorems 1.1. and 1.2. With the
help of the Killing form we identify the Cartan subalgebra $\hog$
with its conjugate space. This identification allows to define the
Killing form $(\cdot,\cdot)$ on $\hog^*$. Note that the  condition
$\sg^2=\id$ implies $(\xi_i,\xi_j) = 0$ for different elements
$\xi_i,\xi_j$ in $\SC$.
\\
{\bf Lemma 2.1}. Let $\eta\in\Dp$ and $t=j(\eta)$. We claim: \\
 1) if  $\sg_{t-1}(\eta) > 0$ (i.e $\eta\in\Pi^{(t)}$), then $\sg_k(\eta)>0$ for all
 $0\lee k \lee t-1$;\\
2) if $\sg_{t}(\eta) > 0$, then $\sg_{t-1}(\eta) > 0$.\\
 {\bf Proof}. We shall prove statement 1). The statement 2) is proved similar.
If  $\eta\in \SC$, then $r_m(\eta)=\eta$  for all $m$. This implies
1).

Suppose that  $\eta\in\Dp\setminus \SC$. Then the only one of the
following cases will take place:
\\
i)~$(\eta, \xi)=0$ for all  $\xi\in \SC$;\\
ii)~ there exists a unique  $\xi'\in \SC$ such that  $(\eta,
\xi')\ne 0$ and  $(\eta, \xi)=0$
for all  $\xi\in \SC\setminus\{\xi'\}$;\\
iii)~ there exist exactly two  $\xi',\xi''\in\SC$ such that $(\eta,
\xi')\ne 0$, $(\eta, \xi'')\ne 0$ and $(\eta, \xi)=0$ for all
$\xi\in \SC\setminus\{\xi',\xi''\}$.

The statement 1) is trivial for the case i). In the case ii) we
embed $y_\eta$ in to the subalgebra
$\mathrm{span}\{y_\gamma\vert~\gamma\in\Dp, \gamma\in
\Zb\eta+\Zb\xi'\}$ that is isomorphic to  $\ut(3,K)$. The proof is
reduces to the  special case $\nog=\ut(3,K)$ and $\SC=\{\alpha\}$
where $\alpha$ is one of positive roots $\alpha_{12}, \alpha_{13},
\alpha_{23}$.

In the case iii) embed  $y_\eta$ in to subalgebra
$\mathrm{span}\{y_\gamma\vert~\gamma\in\Dp, \gamma\in
\Zb\eta+\Zb\xi'+ \Zb\xi''\}$ that is isomorphic to $\ut(4,K)$. The
proof is reduces to the  special case   $\nog=\ut(4,K)$ and
 $\SC$
coincides with one of the following subsets  $\SC_1= \{\alpha_{14},
\alpha_{23}\}$, $\SC_2= \{\alpha_{13}, \alpha_{24}\}$, $\SC_3=
\{\alpha_{12}, \alpha_{34}\}$. $\Box$
\\
 {\bf Lemma 2.2}. For any $\gamma'\in \CC_-$ there exists a unique  $\gamma\in \CC_+$
such that  $\gamma + \gamma'\in \SC$. The correspondence
$\gamma'\mapsto\gamma$
is bijective. \\
 {\bf Proof} is proceeded similar to Lemma 1 reducing to the unitriangular Lie subalgebra on lower size.
 $\Box$ \\
{\bf Proof of the Theorem 1.1}.
\\
1) Let us show that  $\pog_\sg$ is a subalgebra. It suffices to
prove that, if  $\gamma,\gamma'\in\Pi$ and  $\gamma+\gamma'\in\Dp$,
then  $\gamma+\gamma'\in\Pi$. Let $j(\gamma)=t$, $j(\gamma')=t'$ and
$t < t'$. Then  $j(\gamma+\gamma')=t$. Since $\gamma\in\Pi$, then
$\sg_{t-1}(\gamma)>0$. From Lemma 2.1 we conclude
$\sg_{t-1}(\gamma')>0$.
Therefore,  $\sg_{t-1}(\gamma+\gamma')>0$ and $\gamma+\gamma'\in\Pi$.\\
2) Let us show that  $\pog_\sg$  is an isotropic subspace. It
suffices to show that   $(P+P)\cap \SC = \emptyset.$ Assume the
contrary:~  $\xi=\gamma+\gamma'$ where $\xi\in \SC$, $\gamma,
\gamma'\in\Pi$. Let $t=j(\gamma) < t'=j(\gamma')$. Lemma 2.1 implies
$\sg_t(\gamma')
> 0$. On the other hand,  since $\gamma\in\Pi$, then
$$ \sg_t(\gamma') = \sg_{t-1}r_\xi(\gamma') = -\sg_{t-1}(\gamma) <
0.$$
 A contradiction. \\
3) Let us show that  $\pog_\sg$ is maximal isotropic subspace.
Assume the contrary.
 Let  $\pog_\sg+Ka$ be an isotropic subspace where
$$a=\sum_{\gamma'\in \CC_-} c_{\gamma'}y_{\gamma'}\ne
0.$$

Suppose that  $c_{\gamma'_0}\ne 0$ for some $\gamma'_0\in \CC_-$.
There exists a root $\gamma_0\in \CC_+\subset\Pi$ such that
$\gamma_0+\gamma_0'\in \SC$ and $\gamma_0+\gamma'\notin \SC$ for all
$\gamma'\in \CC_-\setminus \{\gamma_0'\}$.
 Then
$$0 = f([y_{\gamma_0}, a]) =
f(\sum_{\gamma'\in \CC_-}
 c_{\gamma'}y_{\gamma_0 +
\gamma'}) =
c_{\gamma'_0}
 f(y_{\gamma_0 + \gamma_0'}).$$
This contradicts to  $f(y_{\gamma_0 + \gamma_0'})\ne 0$. Finally,
the points  1,2,3  conclude that  $\pog_\sg$ --
is a polarization for all  $f\in X_\sg$. $\Box$ \\
{\bf Proof of the Theorem  1.2}. The dimension of any maximal isotropic subspace of the skew symmetric form
 $< x,y>_f =
f([x,y])$ equals to  $\frac{1}{2}(\dim\nog+\dim\nog^f)$ where
$\nog^f$ --
 is a Lie algebra of stabilizer of $f$.
 Hence
$$\dim\Omega(f) = \dim\nog-\dim\nog^f = 2 ~\mathrm{codim}~ \pog_\sg =
2\vert \CC_-\vert = \vert \CC_+\vert + \vert \CC_-\vert =
\vert\Dp\vert - \vert \MC \vert -\vert \SC\vert.\eqno(4)$$

By Remark 4, $\vert \MC\vert = \vert\Dp_\sg\vert$. Therefore, $\vert
\MC\vert = \vert\Dp\vert - l(\sg)$. Substituting in (4) we have got
$\dim\Omega(f)=l(\sg)-\vert\SC\vert = l(\sg)-s(\sg)$.$\Box$

\section{Defining ideal of  orbit}

In this section we shall prove Theorem 1.4 on the defining ideal of
an orbit associated with involution. Let $f\in X_\sg$. Denote by
$\IC$  the ideal in the symmetric algebra $S(\nog)$ generated by
$Q_{i,t}$,  $\eps_t-\eps_i\in\MC$, and $D_m - D_m(f)$,
   $ 1\le m\le s$.

Our goal is to prove that  $\IC$ coincides with  the
 defining ideal  $\IC(\Omega(f))$ of the orbit  $\Omega(f)$. The proof is divided into several steps :\\
1) the ideal  $\IC$ annihilated at  $f$ (Proposition 3.1));\\
 2) the ideal  $\IC$ is prime  and  $\dim(\Ann\IC)=\dim\Omega$. (Proposition 3.7) \\
 3) the ideal $\IC$ in $S(\nog)$ is stable with respect to the adjoint action of the group
 N (Proposition  3.8).\\
 {\bf Proposition 3.1}  The ideal  $\IC$ annihilates at   $f$.\\
 {\bf Proof}. Easy to see that  the elements $D_m-D_m(f)$ annihilate at $f$.
If the square $(i,t)$, $i>t$, is filled by symbol $\bullet$ and has
type 0, then $Q_{i,t}= D_{i,t}$. All entries of $i$th row of
$D_{i,t}$ annihilate at $f$. Hence $Q_{i,t}(f)=0$.

Suppose that  the $(i,t)$,~ $i>t$ is filled by symbol $\bullet$ and
 has type 1. As in
Remark 5  $k=\sg_{t-1}(i)=\sg(i)$,~ $a=\sg(t)$ and $a< k < t < i$.
Put $J=J(k,t)$, $I= I(k,t)$ (see \S1).
   Decompose  $$\begin{array}{c} J=\{a,t\}\sqcup\La_1\sqcup\La_2,~\mbox{where }\\
   \La_1=\{ j\in J:~ k< \sg(j) < t\}\\
$$\La_2 = \{ j\in J:~  \sg(j) > t\}.\\\end{array}\eqno(5)$$
Present $\La_1$ as a union $\La_1 =\La_{1,0}\sqcup \La_{1,1}$, where
$$\La_{1,0}=\{j\in\La_1:~ 1\lee j <k\},$$
$$ \La_{1,1}=
\{j\in\La_1:~ k<j <t\}.$$
 The subset   $\La_{1,1}$ is stable with
respect to action of $\sg$.

Note that since $\sg(i)=t$, then $i=\sg(t)$ and $k\in\La_2$. We have
$$\begin{array}{c} J = J(k,t)= \{a,t\}\sqcup \La_{1,0}\sqcup \La_{1,1}\sqcup \La_2,\\
I = I(k,t)= \{k,t\}\sqcup \sg\La_{1,0}\sqcup \La_{1,1}\sqcup \sg\La_2,\\
I\cap J = \{k,t\}\sqcup \La_{1,1}\end{array}\eqno(6)$$

As in (3) we decompose  $D_{k,t}(\tau)$ in powers of  $\tau$. Denote
$$\begin{array}{c} J_0=J\setminus I = \{a\}\sqcup \La_{1,0}\sqcup \La_2\setminus\{k\},\\
I_0=I\setminus J = \sg\La_{1,0}\sqcup \sg\La_2.\end{array}\eqno(7)$$
We obtain
 $$Q_{i,t}= P_{k,t,1} = \sum_{j\in I\cap J} \ueps(j)
M^{J_0\sqcup\{j\}}_{I_0\sqcup\{j\}}= \ueps(t)
M^{J_0\sqcup\{t\}}_{I_0\sqcup\{t\}} + \sum_{j\in I\cap J,~j\ne t}
\ueps(j) M^{J_0\sqcup\{j\}}_{I_0\sqcup\{j\}},\eqno(8)$$ where
$\ueps(j)=\pm 1$. The  minor $M^{J_0\sqcup\{t\}}_{I_0\sqcup\{t\}}$
annihilates at $f$ since all its  $t$th column annihilates are equal
to zero at $f$.
 Each minor $ M^{J_0\sqcup\{j\}}_{I_0\sqcup\{j\}}$,~ $j\in I\cap J$ and $j\ne t$,
  annihilates at
 $f$ since its $a$th column  annihilates at $f$.
We have got  $Q_{i,t}(f)=0$.
 $\Box$.\\
 {\bf Remark 6}.  If the square $(i,t)$ is filled by symbol $\bullet$ and  has type
 1, then the polynomial  $P_{k,t,0}$ coincides with
 $M_{I_0}^{J_0}$. This minor  annihilates  at  $f$ since
 all entries of its  $a$th column  annihilates at $f$.

Let $A=(a_{ij})$  be a $n\times n$-matrix where entries $a_{ij}$
are elements of some  commutative domain $\CC$.\\
 {\bf Lemma 3.2}.  Let $I,J$ be  two subsystems of
$\{1,\ldots, n\}$ and $|I|=|J|$.
Suppose that $J=J_1\sqcup J_2$, $I = I_1\sqcup I_2$, where $|I_1|=|J_1|$.\\
1) Let $S\subset J_2$. The following identity holds for   minors of
the  matrix $A$:
 $$ M_{I_1}^{J_1}M_I^J = \sum_{T\subset I_2,~\vert T\vert
 =\vert S\vert} \ueps(T,S) M_{I_1\sqcup T }^{J_1\sqcup S }
 M^{J\setminus S}_{I\setminus T}.\eqno(9)$$
2) Let $T\subset I_2$. The following identity holds for   minors of
the  matrix $A$:
 $$ M_{I_1}^{J_1}M_I^J = \sum_{S\subset J_2,~\vert T\vert
 =\vert S\vert} \ueps(T,S)
  M_{I_1\sqcup T }^{J_1\sqcup S } M^{J\setminus S}_{I\setminus T}.\eqno(10)$$
 Here  $\ueps(T,S)=\pm 1$. The sign   depends on the choice of $T, S, I_1, I_2, J_1, J_2$.
 \\
{\bf Proof}. We may assume that $a_{ij}$ are variables. Let $M_I^J$
be the minor of the  matrix $A$ with the system of columns $J$ and
the system of  rows $I$.
 Consider the determinate $$\widetilde{M}_{I_2}^{J_2}=\det\left(\tilde{a}_{ij}\right)$$ of size $|I_2|= |J_2|$ with entries
 $$\tilde{a}_{ij} = M^{J_1\sqcup\{j\}}_{I_1\sqcup\{i\}},~~\mbox{where}~~ j\in J_2, i\in I_2.$$
Apply the Laplace formula:
$$ \widetilde{M}_{I_2}^{J_2} = \sum_{T\subset I_2,~|T|=|S|}\pm\widetilde{M}_{T}^{S}
\widetilde{M}_{I_2\setminus T}^{J_2\setminus S}\eqno(11)$$ The
following formula holds
$$
\left(M_{I_1}^{J_1}\right)^{|J_2|-1}M_{I}^{J} = \pm
\widetilde{M}_{I_2}^{J_2}\eqno(12)$$ Let us prove (12). By
elementary transformations  of $M_I^J$ one can obtain
 zeroes in   $I_1\times J_2$ and $I_2\times J_1$ preserving
  entries of   $I_1\times J_1$.   Then
  $$M_I^J=\pm M_{I_1}^{J_1} \overline{M}_{I_2}^{J_2}$$
  where $ \overline{M}_{I_2}^{J_2} = \det(\overline{a}_{ij})$,
  $\overline{a}_{ij}\in\mathrm{Fract}(\CC)$
  and $\tilde{a}_{ij}=\pm \overline{a}_{ij}M_{I_1}^{J_1}$.
  We have got
  $$ \widetilde{M}_{I_2}^{J_2} = \pm \left(M_{I_1}^{J_1}\right)^{|J_2|}
  \overline{M}_{I_2}^{J_2} =
  \pm \left(M_{I_1}^{J_1}\right)^{|J_2|-1} M_{I_1}^{J_1}\overline{M}_{I_2}^{J_2} =
  \left(M_{I_1}^{J_1}\right)^{|J_2|-1}M_{I}^{J}. $$
  This proves (12).
 Substituting (12) for (11) we have got:

$$
\left(M_{I_1}^{J_1}\right)^{|J_2|-1}M_{I}^{J} =  \sum_{T\subset I_2,~|T|=|S|}\ueps(T,S)
\left(M_{I_1}^{J_1}\right)^{|S|-1}M_{I_1\sqcup T}^{J_1\sqcup S}
\left(M_{I_1}^{J_1}\right)^{|J_2\setminus S|-1}M_{I\setminus T}^{J\setminus S}
$$
Cutting  $\left(M_{I_1}^{J_1}\right)^{|J_2|-2}$ in the left and
right sides of this  equality, we obtain (9). The statement 2) is
proved similarly. $\Box$

Let $C=\tau^m(c_0+c_1\tau+c_2\tau^2+\ldots)$, ~$c_0\ne 0$, be a
polynomial in $\tau$ with coefficients in some commutative domain
$\CC$.
We call  $m$ the lower  degree of $C$. Easy to see that  $m(AB)= m(A)m(B)$.\\
{\bf Definition 3.3}. We call the presentation $C$ as a sum
$C=A_1+\ldots+A_k$ of polynomials $A_i\in\CC[\tau]$  an  admissible
presentation, if
$ m(C)\lee m(A_i)$ for all nonzero   $A_i$.\\
 {\bf Remark 7}.  Let $D= \tau^m(d_0+d_1\tau+d_2\tau^2+\ldots)$,
$A_i=\tau^{m(A_i)}(a_{i,0}+a_{i,1}\tau+\ldots)$,~
$B_i=\tau^{m(B_i)}(b_{i,0}+b_{i,1}\tau+\ldots)$ and $C$ as above.
Let $D, C\ne 0$.  If $DC=\sum A_iB_i$ is an admissible presentation,
then
 $$d_0c_0=\sum'
a_{i,0}b_{i,0},\eqno(13)$$
$$d_0c_1 = - d_1c_0 +  \sum' (a_{i,1}b_{i,0}+
a_{i,0}b_{i,1}) + \sum'' a_{i,0}b_{i,0},\eqno(14)$$ where  the sum
$\sum'$ means the sum over  those $i$ that  obeys
$m(D)+m(C)=m(A_i)+m(B_i)$; the sum  $\sum''$ means the sum over
those $i$
that obeys  $m(D)+ m(C)+1= m(A_i)+m(B_i)$.\\
{\bf Remark 8}. Let $A,B,C,D$ as above. If all $a_{i,0}$,~ $a_{i,1}$
lie in some ideal $\IC$ of $C$ and if $d_0$ is invertible modulo
$\IC$, then $c_0$,~ $c_1 $ lie in $\IC$.\\
{\bf Definition 3.4}. Let $v=(v_1<\ldots < v_n)$ and
$u=(u_1<\ldots<u_n)$ be two increasing $n$-dimensional vectors. We
say $v\lee u$, if
 $v_i\lee u_i$ for all $i$. Respectively, $v<u$, if $v\lee u$ and $v\ne
 u$.
 For  $I, J\subset \{1,\ldots, n\}$, $|I|=|J|$, we say $I\gee J$ if
 $\mathrm{ord}(I)\gee \mathrm{ord}(J)$ in the above sense.\\
 {\bf Lemma 3.5}. Let $I_1$, $J_1$, $I$, $J$ be as in
Lemma 3.2. Suppose that $I\gee J$ and $I_1\gee J_1$. We claim that
the decompositions   (9) and (10) are admissible over $S(\nog)$.\\
{\bf Proof} is given  for (9). The case of formula (10) is similar.
Under the assumptions the left  side of (9) is nonzero (see Remark
1).  Let $m$ (resp. $m_T$) be the lower degree of the left part of
(9) (resp.  of the nonzero $T$th summand   in the right side) Let us
show that $m\lee m_T$.  By definition,
$$m=|J_1|+|J| - |I_1\cap J_1|-|I\cap J|$$
$$m_T = |J_1\sqcup S|+|J\setminus S| - |(I_1\sqcup T)\cap (J_1\sqcup S)|-
|(I\setminus T)\cap (J\setminus S)|.$$
Since
$$|J_1|+|J|=|J_1\sqcup S|+|J\setminus S|,$$
$$|(I_1\sqcup T)\cap (J_1\sqcup S)| = |I_1\cap J_1| + |I_1\cap S|
+|T\cap J_1|+|T\cap S|,$$
$$(I\setminus T)\cap (J\setminus S)| = |I\cap J| - |I\cap S| -|T\cap J|
+ |T\cap S|,$$ then $m_T-m=  |I\cap S| + |T\cap J| - |I_1\cap S| -
|T\cap J_1| - 2|T\cap S| = |I_2\cap S|+|T\cap J_2| - 2|T\cap S|
=|T\cap(J_2\setminus S)| + |(I_2\setminus T)\cap S|\gee 0$. $\Box$\\
 {\bf Lemma 3.6}. Let $J\subset\{1,\ldots, n\}$.
Suppose that $J$ have the following property: if $j\in$ and $j'<j$,
~$\sg(j')>\sg(j)$, then $j'\in J$.
 Denote $P_\mu = P^J_{I,\mu}$,
  $\mu\in\{0, 1\}$ (see (2)). Let $I\subset \{1,\ldots,n\}$ and $|I|=|J|$.
  We claim,
  \\
  1) if $I=\sg(J)$, then $P_0-P_0(f)\in \IC$, $ P_0(f)\ne 0$ for any
  $f\in X_\sg$;\\
  2) if $I>\sg(J)$, then $P_\mu\in \IC$ for $\mu\in\{0, 1\}$.
\\
{\bf Proof}. First, note that the assumption on $J$, imply that
$\sg(J)\gee J$ in the sense of Definition 3.4. By Remark 1,
$M_I^J(\tau)\ne 0$ and  its lower degree equals to $|J\setminus
I|=|I\setminus J|$.

 We prove the statement using the method of induction on
$|I|$. Easy to see that the claim is true for $|I|=1$. Assume that
the claim is true
for $|I|< m$. We a going to prove the claim for $|I|=m$.\\
{\bf Case 1}.~ Let $I=\sg(J)$. Easy to see that $P_0(f)\ne 0$.  Put
$t=\max J$, ~$\sg(t)=a$. Denote
$$
J_1 = \{ j\in J:~ \sg(j) > a \}, \quad I_1=\sg(J_1).$$
  {\bf Point 1a}.~ $|J_1| < |J|-1$.
Denote $$ M^{(1)}(\tau) = M_{I_1}^{J_1}(\tau),\quad
M(\tau)=M_I^J(\tau),\quad M^{(\beta)}(\tau) =
M_{I_1\sqcup\{\beta\}}^{J_1\sqcup\{t\}}(\tau), \quad
M^{(*\beta)}(\tau) = M_{I\setminus\{\beta\}}^{J\setminus \{t\}}
(\tau).$$

If $J_1=\emptyset$, then we put $M^{(1)}(\tau)=1$.
 Write down the admissible decomposition (9)
 with $S=\{t\}$ for  matrix $\Phi(\tau)$:
$$  M^{(1)}(\tau) M(\tau)=
\sum_{\beta\in I(\tau),~ \beta\lee a} \ueps(\beta)
M^{(\beta)}(\tau)M^{(*\beta)}(\tau)=$$
$$ \pm
M^{(a)}(\tau)M^{(*a)}(\tau) + \sum_{\beta\in I(\tau),~ \beta > a}
\ueps(\beta) M^{(\beta)}(\tau)M^{(*\beta)}(\tau).\eqno(15)
$$
Denote by $P^{(1)}_0$, $P^{(\beta)}_0$, $P^{(*\beta)}_0$ the lowest
coefficients of  $M^{(1)}(\tau)$, $M^{(\beta)}(\tau)$,
$M^{(*\beta)}(\tau)$. Recall $P_0$ is the lowest coefficient of
$M(\tau)$.
 By (13), we have got
$$P^{(1)}_0 P_0 =\sum'_{\beta\in I(\tau),~ \beta\lee a} P^{(\beta)}_0
P^{(*\beta)}_0.\eqno(16)$$

The subset $J_1$ also obey assumption of Lemma.  By induction
assumption,  $P^{(1)}_0-P^{(1)}_0(f)\in\IC$, $P^{(1)}_0(f)\ne 0$.
The left side of (16) is nonzero  at $f$.

The subset  $J\setminus\{t\}$ also obey the assumption of Lemma.
 If $\beta<a$,
then $I\setminus\{ a\} < I \setminus\{\beta\}$ (in the sense of
Definition 3.4). By the induction assumption, $P^{(*\beta)}_0\in\IC$
for $\beta < a$  By Lemma 3.1, $P^{(*\beta)}_0(f)=0$. Considering
values  of right and left parts of (16) at $f\in X_\sg$, we conclude
that the term $\beta=a$ occurs in $\sum'$:
$$P^{(1)}_0 P_0 = \pm P^{(a)}_0
P^{(*a)}_0 + \sum''_{\beta\in I, \beta < a}  P^{(\beta)}_0
P^{(*\beta)}_0.\eqno(17)$$ By the induction assumption, $\pm
P^{(a)}_0$,~$\pm P^{(*a)}_0$ satisfy 1). That is $P-P(f)\in\IC$,
$P(f)\ne 0$  for any of this two polynomials. This proves
$P_0-P_0(f)\in \IC$ and $P_0(f)\ne 0$ in Point 1a.\\
{\bf Point 1b.}~$|J_1| = |J|-1$. As above $t=\max J$,~
$a=\sg(t)=\min I$. In Point 1b  $J=J(a,t)$,~ $I=I(a,t)$ (see \S 1).
 If $a>t$,  then $P_0$ is one of the
minors $D_p$, $1\lee p\lee s$. This proves statement 1).

Let $a\lee t$. Define $I_1$ and $J_1$ as in Point 1a. In Point 1b~
$J_1=J-\{t\}$, $I_1=I-\{a\}$ and $P_0 = P^{(1)}_0$. Subsets $I_1$,
$J_1$ obey induction assumption. The statement 1) is true for
$P^{(1)}_0$. Hence 1) is true for $P_0$.\\
{\bf Case 2.}~ $I>\sg(J)$. There exists $k\in I$ and $ b\in J$ such
that $k\notin \sg(J)$ and $k>\sg(b)$. Similar to Case 1, we denote
$$
J_1 = \{ j\in J:~ \sg(j) > k \}, \quad I_1=\sg(J_1), \quad
P^{(1)}=P^{J_1}_{I_1, 0}.$$
 {\bf Point 2a.}~ $|J_1| <
|J|-1$.  Write down the admissible decomposition (10)
 with $T=\{k\}$ for  matrix $\Phi(\tau)$:
$$
M_{I_1}^{J_1}(\tau) M_I^J(\tau)=\sum_{\alpha\in J\setminus J_1}
\ueps(k) M_{I_1\sqcup\{k\}}^{J_1\sqcup\{\alpha\}}(\tau)
M_{I\setminus\{k\}}^{J\setminus \{\alpha\}} (\tau)\eqno(18)$$ The
pairs of subsets $I_1, J_1$ and   $I_1\sqcup\{k\}$,
$J_1\sqcup\{\alpha\}$ obey the conditions of Lemma 1) and 2)
respectively.  Hence, $P^{(1)}_0 - P^{(1)}_0(f)\in\IC$,
$P^{(1)}_0(f)\ne 0$ and
$$P^{J_1\sqcup\{\alpha\}}_{I_1\sqcup\{k\}, \mu}\in\IC$$ for $\mu\in\{0, 1\}$.
By (18) and Remark 8, $P_\mu\in\IC$.\\
{\bf Point 2b}.~ $|J_1| = |J|-1$. In this case $I=I_1\sqcup\{k\}$
and $k=\min I$. As above $t=\max J$.

If $b<t$, then we consider $J_2=\{j\in J:~ 1\lee j < b\}$ and $I_2=
\sg J_2$. We write down (9) for subsets $I_2\subset I$,~$J_2\subset
J$ and $S=\{b\}$. We conclude the proof similar to Point 2a.

Let $b=t$. It follows  $J=J(k,t)$, ~ $I=I(k,t)$. If $k>t$, then the
square $(k,t)$ is filled by  symbol $\bullet$ in the admissible
diagram; the square $(k,t)$ has type 0. The minor $M_I^J(\tau)$ is a
homogeneous polynomial in $\tau$. We have $P_0=Q_{k,t}\in\IC$ and
$P_1=0$. If $k\lee t$, then the square $(i,t)$, where $i=\sg(k)$, is
filled by symbol $\bullet$ in the admissible diagram; the square
$(i,t)$ has type 1.  This proves $P_1 = Q_{i,t}\in\IC$.

To prove that $P_0\in\IC$ we recall that the pair $I_1, J_1$ obeys
the induction assumption. Denote    $(I_*= I_1\setminus \{t\})\sqcup
\{k\}$. Note $I_*> I_1$ and $P_0=P^{J_1}_{I_*,0}$. The induction
assumption concludes the proof. $\Box$
\\
 {\bf  Proposition 3.7}. \\
   $\IC$ is a prime ideal and $\dim(\Ann\IC)=\dim\Omega(f)$.\\
{\bf Proof}.
 Consider the linear order in the set of squares $\{(i,t):~t<i\}$
 of admissible diagram such that
 that\\
 1) if $t_1<t$, then $(i_1,t_1)<(i,t)$;\\
 2) if the square  $(i,t)$ is filled by symbol $\bullet$ or symbol  $\otimes$,
 then
 $(i,t)>(i_1,t)$ for any square $(i_1,t)$  filled by symbol  $+$ or $-$.\\
 3) $(i_1,t)<(i,t)$  if $i_1 > i$ in any other case that do not mentioned
in 1) and 2).

We shall construct a system of elements    $\{\tilde{y}_{it}|~
i>t\}$ such
that\\
i) any $\tilde{y}_{it}$ has the form $\tilde{y}_{it} =
y_{it}+v_{it}$ where $v_{it}$ lies in subalgebra in $S(\nog)$
generated by $y_{i_1,t_1}$, $(i_1,t_1)<(i,t)$; \\
ii)  the ideal $\IC$ is generated by  elements
 $\tilde{y}_{it} -
 \tilde{y}_{it}(f)$, where  the square $(i,t)$ filled
 by symbol $\otimes$, and $\tilde{y}_{it}$, where he square $(i,t)$ is filled by
 symbol  $\bullet$.

The condition i) implies that the system of elements
$\{\tilde{y}_{it}|~ i>t\}$ generate the symmetric algebra $S(\nog)$
and ii) implies the factor algebra of $S(\nog)$ over the ideal $\IC$
is isomorphic to the algebra of polynomials  generated by
$\{y_{it}\}$, where  the square $(i,t)$ is filled by $+$ or $-$ in
the admissible diagram. This proves that the ideal $\IC$ is prime
and $\dim(\Ann\IC)$ is equal to the number of $\pm$'s in the
admissible diagram. By Remark 4, $\dim(\Ann\IC)=\dim\Omega(f)$.

So it suffices  to construct the required  system of elements
$\{\tilde{y}_{it}|~ i>t\}$ that obey i) and ii). First, we put
$\tilde{y}_{it}= y_{it}$ if the square $(i,j)$ is filled by symbol
$+$ or $-$.

Suppose that the square  $(i,t)$ is  filled by symbol $\otimes$ or
symbol $\bullet$ of type 0. Decomposing $D_{i,t}$  by the column
$t$, we obtain
 $ D_{k,t}  = Ay_{it}+u_{it}$,  where $u_{it}$ as i) and
    $A$ equals to some nonzero constant $a$ modulo $ \bmod \IC$ (see Lemma 3.6(1)).
 Put $\tilde{y}_{it}= y_{it} +a^{-1}u_{it}$.

Suppose that the square $(i,t)$  is filled by symbol $\bullet$  and
has type 1. As in Remark 5 we denote $k = \sg (i) <t$.
   We preserve the notations (5). Using (6) and (7) we have got (8).
 Hence, $$Q_{i,t} = \pm M^{J_0\sqcup\{t\}}_{I_0\sqcup\{t\}} + u_{it}\eqno(19)$$
 where $u_{it}$ lies in the subalgebra generated by $y_{cd}$,~ $d<t$.
 Let as in Lemma 3.1  $a=\sg(t)$.
Recall  $J_0=\{a\}\sqcup \La_{10}\sqcup(\La_2\setminus k)$, ~$I_0 =
\sg \La_{10}\sqcup\sg\La_2$.
 Decompose
 $\La_2=\La_{2,0}\sqcup\La_{2,1}$, where
 $$
 \La_{2,0} = \{ j\in \La_{2}:~ 1\le j< a\},\quad  \La_{2,1} =
 \{ j\in \La_{2}:~ a< j< t\}.$$
Let $\La_{2,1}=\{k_1,\ldots, k_p\}$ and $\sg \La_{2,1}=\{i_1,\ldots,
i_p\}$ where each $i_b=\sg(k_b)$. By definition, $k\in \La_{2,1}$
and $i\in \sg\La_{2,1}$. Note that all squares $\{i_m,t\}$ are
filled by symbol $\bullet$ and have type 1.
 For
convenience of reader we put just  here \\
{\bf Example 4.} The diagram below is the  $10\times 10$ admissible
diagram of  involution $\sg=(1,10)(2,5)(3,7)(4,9)(6,8)$. We put
number of columns in the squares of first row and the numbers of of
rows in the first column. Here we write decompositions for $t=7$,~
$i=9$. Hence, $a=3$  and $k=4$.

\begin{center} {
\begin{tabular}{|p{0.25cm}|p{0.25cm}|p{0.25cm}|p{0.25cm}|p{0.25cm}|p{0.25cm}|
p{0.25 cm}|p{0.25cm}| p{0.25cm}|p{0.25cm}|p{0.25cm}|}
\hline &1 &2 &3 &4 &5 &6 &7&8&9&10 \\
\hline 1&    &           &           & & & &&&&\\
\hline 2&$+$ &           &           & & & &&&&\\
\hline 3&$+$ & $+$       &           & & &&&&&\\
\hline 4&$+$ & $+$       & $+$       &  & & &&&&\\
\hline 5&$+$ & $\otimes$ & $-$       & $-$ &          &    & &&&\\
\hline 6&$+$ & $\bullet$ & $+$       & $+$ & $\bullet$     &          & &&&\\
\hline 7&$+$ & $\bullet$ & $\otimes$ & $-$ & $\bullet$     & $-$      & &&&\\
\hline 8&$+$ & $\bullet$ & $\bullet$ & $+$ & $\bullet$     & $\otimes$& $\bullet$&&&\\
\hline 9&$+$ & $\bullet$ & $\bullet$ & $\otimes$ & $\bullet$& $-$     &$\bullet$ &$-$&&\\
\hline 10&$\otimes$ & $-$ & $-$ & $-$ & $-$& $-$                      &$-$&$-$&$-$&\\
\hline
\end{tabular}
\quad\quad $\begin{array}{l}\La_1 =\La_{10} = \{2\},~~\sg \La_{10},
=\{5\}\\
\\
 \La_{20}=\{1\},~~\sg \La_{20}= \{10\},\\
\\
  \La_{21}=\{6\},~~\sg \La_{20}= \{8\},\\
  \\
  J_0=\{1,2,3, 7\}, ~~ I_0=\{5,8,9,10\},\\
\\
M^{J_0}_{I_0} = M^{1,2,3,6}_{5,8,9,10}.\\
\\
Q_{9,7} = -
 M^{1,2,3,6,7}_{5,7,8,9,10} - M^{1,2,3,4,6}_{4,5,8,9,10}
\end{array}$
}
\end{center}
The last formula for $Q_{9,7}$ is the decomposition (19) for $i=9$
and $t=7$. End of Example.

We continue the proof. Applying the formula (9) with
$S=\La_{1,0}\sqcup \{a\}$ for the matrix $\Phi$  we obtain
 $$
 M^{\La_{2,0}}_{\sg\La_{2,0}} M^{J_0\sqcup\{t\}}_{I_0\sqcup\{t\}} =
 \sum_{T\subset \sg\La_{1,0}\sqcup \sg\La_{2,1}\sqcup\{t\}} \pm
 M^{\La_{2,0}\sqcup \La_{1,0}\sqcup\{a\} }_{\sg\La_{2,0}\sqcup T}
 M^{(\La_2\setminus \{k\})\sqcup\{t\}}_{I_0\sqcup\{t\}\setminus T}\eqno(20)$$
 Since $t$ is greater than any number of  $\sg\La_{1,0}$ and smaller than
 any  number of $\sg\La_{2,1}$, then
all minors $$M_T := M^{\La_{2,0}\sqcup \La_{1,0}\sqcup\{a\}
}_{\sg\La_{2,0}\sqcup T}$$ for $T\subset \sg\La_{1,0}\sqcup
\sg\La_{2,1}\sqcup\{t\}$ and $T\ne\sg\La_{1,0}\sqcup\{t\}$ lie in
the ideal  $\IC$. We obtain
 $$
 M^{\La_{2,0}}_{\sg\La_{2,0}} M^{J_0\sqcup\{t\}}_{I_0\sqcup\{t\}} =
 \pm M^{\La_{2,0}\sqcup \La_{1,0}\sqcup\{a\} }_{\sg\La_{2,0}\sqcup \sg\La_{1,0}\sqcup\{t\}}
 M^{(\La_2\setminus \{k\})\sqcup\{t\}}_{\sg\La_2}\bmod\IC.\eqno(21)
$$

The first minors in left and right sides in (21) are nonzero
constants modulo the ideal $\IC$ (see Lemma 3.6(1)). We have
$$M^{J_0\sqcup\{t\}}_{I_0\sqcup\{t\}} =
\mathrm{const}\cdot M_i,~~\mbox{~where~} M_i = M^{(\La_2\setminus
\{k\})\sqcup\{t\}}_{\sg\La_2}\bmod\IC.\eqno(22)$$
 Substituting (22) for (19) we have got
 $$Q_{it}=c_{i,t}M_i +u_{i,t}\bmod\IC ~\mathrm{~where~} c_{i,t}\in K^*.$$
Then $$c^{-1}_{i,t}Q_{it}=M_i + c^{-1}_{i,t}u_{i,t}.\eqno(23)$$
Decomposing $M_i$ by last $t$th column we obtain
$$M_i =\sum_{r=1}^p
A^{(r)}_iy_{i_r,t} + v_{i,t}\eqno(24)$$
 where $v_{it}$ lies in subalgebra
generated by $y_{c,d}$, ~$d<t$, and $y_{c,t}$, where the square
$(c,t)$ is filled by symbol $-$. Therefore, $v_{i,t}$ lies in
subalgebra generated by $y_{c,d}$, where $(c,d)$ is smaller in the
sense of order defined in the beginning of this proof.

 Note that the matrix
$A=\det(A^{(r)}_{i_m})_{r,m=1}^p$ is invertible modulo ideal $\IC$
(see Lemma 3.6 and formula (12)). Denote $B=(B^{(r)}_m)$ the inverse
matrix for $A$. After, we write down (23) and (24) for $i=i_m$,
~$1\lee m\lee p$, substitute (24) in (23). Further,  multiply (23)
by $B^{(q)}_s$ and summarize:
$$\sum_{q=1}^pB^{(q)}_m c_{i_q,t}^{-1} Q_{i_q,t} =
y_{i_m,t}+w_{i_m,t}\in\IC.$$ Finally, put $\tilde{y}_{i,t}=
y_{i,t}+w_{i,t}$ for any $i\in \sg\La_{21}=\{i_1,\ldots. i_p\}$. The
system $\{\tilde{y}_{i,t}\}$ obey i), ii). $\Box$
\\
 {\bf Proposition 3.8} The ideal  $\IC$  in $S(\nog)$ is stable
 with respect to the adjoint action of group $N$.
 \\
 {\bf Proof}.  The statement of Proposition if equivalent to claim
 that $\IC$ is a Poisson ideal in $S(\nog)$ (with respect to natural Poisson bracket).
 Let $y=y_{p+1,p}$ where $1\lee p\lee n-1$.
 It suffices to prove that $\{y,P\}\in\IC$ for any standard generator $P$.
 Recall that, by definition, $\IC$ is generated by $Q_{i,t}$,  $\eps_t-\eps_i \in \MC$,
  and
$D_m - D_m(f)$,
   $ 1\le m\le s$. \\
 {\bf Case 1}.~ $P$ is equal to $D_m$ or $Q_{i,t}$, where square
 $(i,t)$ is filled by $\bullet$ and has type 0. Then $P$ is a minor
 $ M_I^J$ of  the matrix $\Phi$;~ squares $I\times J$ are  placed under the
 diagonal and $I\gee\sg(J)$.
 Easy to see that
 $$\{y,M_I^J\}= \ueps_1 M_{I_*}^{J} - \ueps_2  M_{I}^{J_*},~~ \mathrm{where}$$
$$\ueps_1= \left\{\begin{array}{l}
1,~~ \mathrm{~if~} p\in I,~p+1\notin I\\
0, ~~\mathrm{otherwise}.
\end{array}\right., \quad \ueps_2= \left\{\begin{array}{l}
1,~~ \mathrm{~if~} p+1\in J,~p\notin J\\
0, ~~\mathrm{otherwise}.
\end{array}\right.,\quad
$$
$$
\begin{array}{l} I_*=(I\setminus \{p\}) \cup\{p+1\}\\
J_*=(J\setminus \{p+1\}) \cup\{p\}
\end{array}.
$$
One can see that $I_*>I\gee \sg(J)$  and $I>\sg(J_*)$. By Lemma 3.2,
minors $M_{I_*}^{J}$,~$M_{I}^{J_*}$ lie in $\IC$. Hence, $\{y,
M_I^J\}\in\IC$.\\
{\bf Case 2}. Suppose that $P=Q_{it}$ and the square $(i,t)$ is
filled
by $\bullet$ and has type 1. \\
{\bf Point 2a}. Denote $[t]=\{1,\ldots,t\}$ and $P_\mu =
P_{I,\mu}^{[t]}$ where $\mu\in\{0, 1\}$.   By Lemma 3.6, $P_\mu$, ~
$\mu\in\{0, 1\}$, lies in the ideal $\IC$. We shall show in this
point that
$$\{y,P_\mu\}\in\IC, ~~\mu\in\{0, 1\}.$$ More precisely,
$$\{y,P_0\}= c_0 (P_*)_0,\eqno(25)$$
$$\{y, P_1\}=c_1 (P_*)_1 + c_2 y P_0,\eqno(26)$$
where $c_0,c_1,c_2\in\{0, 1, -1\}$,~ $P_* = P^{[t]}_{I_*}$,~ $I_*=
(I\setminus \{p\})\cup\{p+1\}$. In (25) and (26) $c_0=c_1=0$ if the
condition $p\in I$,~ $p+1\notin I$ is not true.

The proof is proceeded separately   in each of 12 cases, concerning
entry of $p$, $p+1$ in $I$, $[t]$ (note  that  the case
$p\notin[t]$, $p+1\in[t]$ is not possible). In other cases proof is
similar.

Here we present the complete proof of (25) and (26) in the case
$p\in I$,~$p+1\notin I$,~$p,p+1\in [t]$.

Suppose that $I$ is an ordered increasing subset. Denote $\La =
[t]$,~$\La_0=\La\setminus I$,~ $I_0 = I\setminus \La$.

The polynomial $P_0$ coincides with  the minor $M^{\La_0}_{I_0}$ .
Since $p, p+1\notin I_0$ and $p+1\in \La_0$, ~ $p\notin \La_0$, then

$$\{y, P_0\}= - M_{I_0}^{(\La_0)_*}, ~~\mathrm{where}~
(\La_0)_*=(\La_0\setminus\{p+1\})\cup\{p\}.$$ The minor
$M_{I_0}^{(\La_0)_*}$ coincides with $(P_*)_0$.

Denote $F_0=\La\cap I$. By assumptions,~ $p\in F_0$, ~$p+1\notin
F_0$.  We have
$$
P_1=\sum_{\alpha\in F_0}\ueps(\alpha) M^{\La_0\sqcup
\{\alpha\}}_{I_0\sqcup\{\alpha\}} = \ueps(p) M^{\La_0\sqcup
\{p\}}_{I_0\sqcup\{p\}} + \sum_{\alpha\in F_0,~\alpha\ne
p}\ueps(\alpha) M^{\La_0\sqcup \{\alpha\}}_{I_0\sqcup\{\alpha\}}.
$$
Here $\ueps(\alpha)=(-1)^{i+j}$ where $i$ equals to the index of
$\alpha$  in $\La$ and $j$ equals to the index of $\alpha$  in $I$.
By direct calculations  we have got
$$\{y, P_1\} = \ueps(p) \left(M^{\La_0\sqcup
\{p\}}_{I_0\sqcup\{p+1\}}\pm y M_{I_0}^{J_0}\right) -
\sum_{\alpha\in F_0,~\alpha\ne p}\ueps(\alpha) M^{(\La_0
\setminus\{p+1\})\sqcup \{p,
\alpha\}}_{I_0\sqcup\{\alpha\}}.\eqno(27)
$$
On the other hand,
$$ (P_*)_1 = \ueps'(p+1) M^{\La_0\sqcup
\{p\}}_{I_0\sqcup\{p+1\}} + \sum_{\alpha\in F_0,~\alpha\ne
p}\ueps'(\alpha) M^{(\La_0 \setminus\{p+1\})\sqcup \{p,
\alpha\}}_{I_0\sqcup\{\alpha\}}.\eqno(28)
$$
 Here $\ueps'(\alpha)=(-1)^{i'+j'}$ where $i'$ equals to the index of
$\alpha$  in $\La$ and $j'$ equals to the index of $\alpha$ in $I_*=
(I\setminus\{p\})\cup\{p+1\}$. Similarly,
 $\ueps'(p+1)=(-1)^{i'+j'}$ where $i'$ equals to the index of
$p+1$  in $\La$ and $j'$ equals to the index of $p+1$  in $I_*$.
Comparing (27) and (28), we conclude that  $\ueps(p)= - \ueps'(p+1)$
and $\ueps(\alpha)=\ueps'(\alpha)$ for $\alpha\ne p$, $\alpha\in
F_0$. We obtain $$\{y, P_1\} = -(P_*)_1\pm y
P_0.$$ This proves (26).\\
 {\bf Point 2b.} We finish the  proof of Case 2.
 Denote as above $k=\sg(i)$. Denote  $J_1=J'(k,t)$, ~$I_1=I'(k,t)$,~ $J(k,t) =
J_1\sqcup\{t\}$, ~ $I(k,t)= I_1\sqcup \{k\}$( see \S 1). Recall that
$Q_{i,t}$ is  the second coefficient of minor with columns $J(k,t)$
and rows $I(k,t)$ of the matrix  $\Phi(\tau)$. That is
$$Q_{i,t}= P^{J_1\sqcup\{t\}}_{I_1\sqcup\{k\}, 1}.$$

  Apply (9) for $J_1\subset [t]$,~ $I_1\subset
 I=\sg[t]\sqcup\{k\}$ and $S= \{t\}$:

 $$
 M_{I_1}^{J_1}(\tau)M_{I}^{[t]}(\tau)=\sum_{\beta\gee k,~\beta\in I}
 \ueps(\beta) M_{I_1\sqcup\{\beta\}}^{J_1\sqcup\{t\}}(\tau)
 M_{I\setminus\{\beta\}}^{[t-1]}(\tau) =$$
 $$ \ueps(k)M_{I_1\sqcup\{k\}}^{J_1\sqcup\{t\}}(\tau)
M_{\sg[t-1]}^{[t-1]}(\tau) + \sum_{\beta>k,~\beta\in I}\ueps(\beta)
M_{I_1\sqcup\{\beta\}}^{J_1\sqcup\{t\}}(\tau)
 M_{I\setminus\{\beta\}}^{[t-1]}(\tau).
  \eqno(29)$$
Denote $m_1=|I_1|=|J_1|$,~$m=|I|=|J|$,~$q_1=|I_1\cap J_1|$,
~$q=|I\cap  J|$.

 According Lemma 3.5, decomposition (29) is an admissible
 decomposition. This means that the lower degree of each summand in
 the right side of (29) is greater than the lower degree of the left
 side.

  We show that the lower degree of the first term of the
  right side coincides with the lower degree of the left side.
  Indeed
  $$\mathrm{ldeg}\left( M_{I_1}^{J_1}(\tau)M_{I}^{[t]}(\tau)  \right)
  =\mathrm{ldeg}\left( M_{I_1}^{J_1}(\tau)  \right) +
  \mathrm{ldeg}\left( M_{I}^{[t]}(\tau)  \right) = (m_1-q_1)+(m-q),$$
 $$\mathrm{ldeg}\left( M_{I_1\sqcup\{k\}}^{J_1\sqcup\{t\}}(\tau) M_{\sg[t-1]}^{[t-1]}(\tau)  \right)
  =\mathrm{ldeg}\left( M_{I_1\sqcup\{k\}}^{J_1\sqcup\{t\}}(\tau)  \right) +
  \mathrm{ldeg}\left( M_{\sg[t-1]}^{[t-1]}(\tau) \right) = $$
  $$[(m_1+1)-(q_1+2)]+[(m-1)-(q-2)]=
  (m_1-q_1)+(m-q).$$

This prove coincidence of lower degree. By Lemma 3.6, the first
coefficients of $M_{I_1}^{J_1}(\tau)$ and $
M_{\sg[t-1]}^{[t-1]}(\tau)$ are invertible modulo $\IC$. Comparing
the coefficients of (29) we have got formulas of type (13)  and
(14). Formula (14) provides an expression of $Q_{i,t}$ as a sum of
polynomials of type $P_{I,\mu}^{[t]}$ and $P_{I',\mu}^{[t-1]}$,
where $I>\sg[t]$ and $I'>\sg[t-1]$, with coefficients in $S(\nog)$.
Applying Point 1a), we conclude $\{y, Q_{i,t}\}\in\IC$. $\Box$

 {\bf The proof of Theorem  1.4} follows directly from Propositions 3.1, 3.7 and
 3.8.


\begin{thebibliography}{100}

\bibitem{K-Orb}
Kirillov A.A., Lectures on the orbit method. Novosibisk, 2002.
\bibitem{K-Var}
Kirillov A.A., Variations on the Triangular Theme,
Amer.Math.Soc.Trans.(2). 1995. V.169. P. 43-72.
\bibitem{K-62}
Kirillov A.A., Unitary representations of nilpotent Lie groups,
Uspekhi Math. Nauk. 1962. V.17, P. 57-110 (russian)
\bibitem{IP}
Ignatev M.V., Panov A.N. Coadjoint orbits for the  group
$\UT(7,K)$, Math.RT/0603649
\bibitem{Mel}
Melnikov A. $B$-Orbits in solutions to the equation $X^2=0$ in
triangular matrices, J.Algebra, V.223,2000,  P.101-108.
\bibitem{D}
Dixmier J. Algebres enveloppantes, Gauthier-Villars, 1974.
\bibitem{GG}
Goto M, Grosshans F., Semisimple Lie algebras, Lectires Notes in
Pure and Applied Mathematics, V.38, Marcel Dekker, INC., New York
and Basel, 1978

\end{thebibliography}
\end{document}